\documentclass[10pt, a4paper]{article}

\usepackage{amsfonts, latexsym, amsmath, amssymb, fancyhdr, amsthm, graphicx, subfigure}
\usepackage[english]{babel}

\newtheorem{theorem}{Theorem}[section]
\newtheorem{lemma}{Lemma}[section]
\newtheorem{proposition}{Proposition}[section]

\newtheorem{remark}{Remark}

\title{\textbf{Modelling hematopoiesis mediated by growth factors: Delay equations describing
periodic hematological diseases}\thanks{Published in Bulletin of Mathematical Biology, vol 68, iss
8, 2321-2351 (2006)}}
\author{Mostafa Adimy$^{a}$, \quad Fabien Crauste$^{a}$ \quad and \quad Shigui Ruan$^{b}$\footnote{\small Research was partially supported by the NSF and the University of Miami.}}
\date{Year 2005}

\begin{document}

\maketitle

\begin{center}
{\large $^{a}$}\emph{Laboratoire de Math\'ematiques Appliqu\'ees, UMR 5142,} \\
\emph{Universit\'e de Pau et des Pays de l'Adour,}\\
\emph{Avenue de l'universit\'e, 64000 Pau, France.}\\
\emph{ANUBIS project, INRIA--Futurs}\\
\emph{E-mail: mostafa.adimy@univ-pau.fr, fabien.crauste@univ-pau.fr}\\
\quad\\
{\large $^b$}\emph{Department of Mathematics, University of Miami,}\\\emph{P. O. Box 249085, Coral
Gables, FL 33124-4250, USA.}\\\emph{ E-mail: ruan@math.miami.edu}
\end{center}

\quad

\begin{abstract}
Hematopoiesis is a complex biological process that leads to the production and regulation of blood
cells. It is based upon differentiation of stem cells under the action of growth factors. A
mathematical approach of this process is proposed to carry out explanation on some blood diseases,
characterized by oscillations in circulating blood cells. A system of three differential equations
with delay, corresponding to the cell cycle duration, is analyzed. The existence of a Hopf
bifurcation for a positive steady-state is obtained through the study of an exponential polynomial
characteristic equation with delay-dependent coefficients. Numerical simulations show that long
period oscillations can be obtained in this model, corresponding to a destabilization of the
feedback regulation between blood cells and growth factors. This stresses the localization of
periodic hematological diseases in the feedback loop.
\end{abstract}

\bigskip{}


\section{Introduction}

Hematopoiesis is the process by which erythrocytes (red blood cells), leukocytes (white blood
cells) and thrombocytes (platelets) are produced and regulated. These cells perform a variety of
vital functions such as transporting oxygen, repairing lesions, and fighting infections. Therefore,
the body must carefully regulate their production. For example, there are $3.5\times 10^{11}$
erythrocytes for each kilogram of body weight, so almost $7\%$ of the body mass is red blood cells.
The turnover rate is about $3\times 10^{9}$ erythrocytes/kg of body weight each day, which must be
carefully regulated by several $O_2$ sensitive receptors and a collection of growth factors and
hormones. Although understanding of blood production process evolves constantly, the main outlines
are clear.

Blood cells, that can be observed in blood vessels, are originated from a pool of hematopoietic
pluripotent stem cells, located in the bone marrow of most of human bones. Hematopoietic
pluripotent stem cells, which are undifferentiated cells with a high self-renewal and
differentiation capacity, give rise to committed stem cells, which form bands of cells called
colony forming units (CFU). These committed stem cells are specialized in the sense that they can
only produce one of the three blood cell types: red blood cells, white cells or platelets. Colony
forming units differentiate in precursor cells, which are not stem cells anymore, because they have
lost their self-renewal capacity. These cells eventually give birth to mature blood cells which
enter the bloodstream.

One can see that the hematopoiesis process is formed by a succession of complex differentiations
from hematopoietic pluripotent stem cells to precursors. These different differentiations,
occurring in the bone marrow, are mainly mediated by growth factors. They are proteins acting, in
some way, like hormones playing an activator/inhibitor role. Each type of blood cell is the result
of specific growth factors acting at a specific moment during the hematopoiesis process.

The red blood cells production, for example, called erythropoiesis, is mainly mediated by
erythropoietin (EPO), a growth factor produced at 90\% by the kidneys. Erythropoietin is released
in the bloodstream due to tissue hypoxia. It stimulates the erythropoiesis in the bone marrow,
causing an increase in circulating red blood cells, and consequently an increase in the tissue
p$O_2$ levels. Then the release of erythropoietin decreases and a regulation of the process is
observed: there is a feedback control from the blood to the bone marrow. In extreme situations,
like bleeding or moving to high altitudes, where needs in oxygen are important, erythropoiesis is
accelerated.

White blood cells are produced during leukopoiesis and the main growth factors acting on their
regulation are Granulocyte-CSF (Colony Stimulating Factor), Macrophage-CSF,
Granulocyte-Macrophage-CSF, and different interleukins (IL-1, IL-2, IL-6, IL-8, etc.). Platelets
are mainly regulated by thrombopoietin (TPO), which acts similarly to erythropoietin.

The hematopoiesis process sometimes exhibits abnormalities in blood cells production, causing the
so-called dynamical hematological diseases. Most of these diseases seem to be due to a
destabilization of the pluripotent hematopoietic stem cell compartment caused by the action of one
or more growth factors. For erythropoiesis, abnormalities in the feedback between erythropoietin
and the bone marrow production are suspected to cause periodic hematological disorders, such as
autoimmune hemolytic anemia. Cyclic neutropenia, one of the most intensively studied periodic
hematological diseases characterized by a fall of neutrophils (white blood cells) counts every
three weeks, is now known to be due to a destabilization of the apoptotic (mortality) rate during
the proliferating phase of the cell cycle.

Mathematical models of hematopoiesis have been intensively studied since the end of the 1970s. To
our knowledge, Mackey \cite{m1978,m1979} proposed the first model of hematopoiesis in 1978 and
1979. This model takes the form of a delay differential equation. Since then it has been modified
and studied by many authors, including Mackey. The works of Mackey and Rudnicki
\cite{mr1994,mr1999} and Mackey and Rey \cite{mr1993, mr1995_2, mr1995_1} deal with age-maturity
structured models of hematopoiesis based on the model of Mackey \cite{m1978}. Recently,
Pujo-Menjouet {\it et al.} \cite{pm2004, pbm} have studied the model of Mackey \cite{m1978,m1979}
and obtained the existence of periodic solutions, with long periods comparing to the cell cycle
duration, describing phenomena observed during some periodic hematological diseases. However, the
influence of growth factors has never been explicitly incorporated in these models.

In 1995 and 1998, Belair {\it et al.} \cite{bmm1995} and Mahaffy {\it et al.} \cite{mbm1998}
considered a mathematical model for erythropoiesis. The model is a system of age and maturity
structured equations, which can be reduced to a system of delay differential equations. They showed
that their model fitted well with experimental observations in normal erythropoiesis but they
stressed some difficulties to reproduce pathological behavior observed for periodic hematological
diseases.

In this paper we consider a system of differential equations modelling the evolution of
hematopoietic stem cells in the bone marrow, of mature blood cells in the bloodstream and of the
concentration of some growth factors acting on the stem cell population. A delay naturally appears
in the model, describing the cell cycle duration. This approach is based on the early work of
Mackey \cite{m1978, m1979} and the recent work of Belait {\it et al.} \cite{bmm1995} and Mahaffy
{\it et al.}  \cite{mbm1998} dealing with erythropoiesis. Our aim is to show that oscillations in
such models are mainly due to the destabilization of the feedback loop between blood cells and
growth factors, causing periodic hematological diseases.

The paper is organized as follows. We first describe the biological background leading to the
mathematical model. After showing the existence of a positive equilibrium, we analyze its local
asymptotic stability. This analysis is performed through the study of a characteristic equation
which takes the form of a third degree exponential polynomial with delay-dependent coefficients.
Using the approach of Beretta and Kuang \cite{bk2002}, we show that the positive steady-state can
be destabilized through a Hopf bifurcation and stability switches can occur. We illustrate our
results with numerical simulations and show that long period oscillations can be observed in this
model, as it can be observed in patients with some periodic hematological diseases.

\section{The Model}

In the bone marrow, hematopoietic stem cells are divided in two groups: proliferating and
non-proliferating, or quiescent, cells. The existence of a quiescent phase, also called $G_0$
phase, in the cell cycle has been proved, for example, by Burns and Tannock \cite{bt}. Quiescent
cells represent the major part of hematopoietic stem cells, that is about 90\% of the hematopoietic
stem cell population. Proliferating cells are cells actually in cycle: they are committed to divide
during mitosis after, in particularly, having synthesized DNA. Immediately after division,
proliferating cells enter the $G_0$ phase where they can stay their entire life.

We denote by $Q(t)$ and $P(t)$ the quiescent and proliferating cell populations at time $t,$
respectively. In the proliferating phase, apoptosis, which is a programmed cell death, controls the
cell population and eliminates deficient cells. We assume that the apoptosis rate, denoted by
$\gamma$, is constant and nonnegative. In the $G_0$ phase, cells can disappear, by natural death,
with a rate $\delta$. They also differentiate in mature blood cells with a rate $g(Q)/Q$, where the
function $g$ is assumed to be nonnegative with $g(0)=0$, because no cell can become mature when
there is no hematopoietic stem cell, and continuously differentiable. Moreover, we assume that the
function $Q\mapsto g(Q)/Q$ is nondecreasing for $Q\geq0$, which is equivalent to
\begin{equation}\label{propg}
0\leq g^{\prime}(0)\leq \frac{g(Q)}{Q}\leq g^{\prime}(Q) \qquad \textrm{ for } Q>0.
\end{equation}
It follows, in particularly, that $g$ is nondecreasing and $\lim_{Q\to+\infty}g(Q)=+\infty$.

Quiescent cells can also be introduced in the proliferating phase during their life in order to
ensure the population renewal, at a nonconstant rate $\beta$. It is generally accepted that $\beta$
depends on the total population of non-proliferating cells \cite{m1978, sachs}. However, the
production of mature blood cells is also mediated by growth factors through the stem cell
population: growth factors induce differentiation and maturation of hematopoietic cells via the
stem cell compartment. Thus, the dependence of $\beta$ on growth factors must be represented. We
assume that $\beta$ is continuously differentiable. Moreover, in the particular case of
erythropoiesis, $\beta$ is an increasing function of the erythropoietin concentration, because a
release of erythropoietin increases the production of red blood cells. Hence, we assume that
$\beta$ is an increasing function of the growth factor concentration, with $\beta(Q,0)=0$, and a
nonincreasing function of the $G_0$ phase population \cite{m1978}.

Thus, the equations modelling the differentiation of hematopoietic stem cells in the bone marrow
are
\begin{eqnarray}
\displaystyle\frac{dQ}{dt}&=&-\delta Q(t) -g(Q(t)) -\beta(Q(t),E(t))Q(t) \vspace{1ex}\nonumber \\
&&+2e^{-\gamma\tau}\beta(Q(t-\tau),E(t-\tau))Q(t-\tau),\vspace{1ex}\label{Q}\\
\displaystyle\frac{dP}{dt}&=&-\gamma P(t) +\beta(Q(t),E(t))Q(t) \vspace{1ex}\nonumber \\
&&-e^{-\gamma\tau}\beta(Q(t-\tau),E(t-\tau))Q(t-\tau),\label{P}
\end{eqnarray}
where $E(t)$ is the growth factor concentration at time $t$. The parameter $\tau>0$ denotes the
average time needed by a proliferating cell to divide, that is, $\tau$ is an average cell cycle
duration. The last term in equation (\ref{Q}) represents the amount of cells coming from the
proliferating phase at division. They are in fact quiescent cells introduced in the proliferating
phase one generation earlier. The factor 2 represents the division of each proliferating cell in
two daughter cells.

At the end of their development, precursors give birth to mature blood cells, which are introduced
in the bloodstream. We denote by $M(t)$ the population of circulating mature blood cells. These
cells only proceed from $G_0$ cells at the rate $g(Q)$. Mature blood cells are degraded, in the
bloodstream, at a rate $\mu\geq0$. Red blood cells usually live an average of 120 days, whereas
platelets live about one week and white blood cells only few hours. Mature blood cell population
satisfies the differential equation
\begin{displaymath}
\frac{dM}{dt}=-\mu M(t) + g(Q(t)).
\end{displaymath}

The growth factor concentration is governed by a differential equation with an explicit negative
feedback. This feedback describes the control of the bone marrow production on the growth factor
production, explained in the previous section. This control acts by the mean of circulating blood
cells: the more circulating blood cells the less growth factor produced. We denote by $f$ the
feedback control. The function $f$ depends on the population of circulating cells $M$ and is
positive, decreasing and continuously differentiable. Then
\begin{displaymath}
\frac{dE}{dt}=-k E(t) +f(M(t)),
\end{displaymath}
where $k\geq0$ is the disappearance rate of the growth factor. In fact, the action of the mature
blood cell population on the production of growth factor is not immediate: it is slightly delayed,
but this delay is negligible in front of the cell cycle duration, so we do not mention it.

At this point, one can notice that system (\ref{Q})--(\ref{P}) is not coupled: the population in
the proliferating phase is not needed in the description of the hematopoiesis process, since
circulating blood cells only come from quiescent cells. From a mathematical point of view, the
solution of (\ref{Q}) does not depend on the solution of (\ref{P}) whereas the converse is not
true. Consequently, we concentrate on the system of delay differential equations
\begin{equation}\label{system}
\setlength\arraycolsep{2pt} \left\{\begin{array}{rcl}
\displaystyle\frac{dQ}{dt}&=&-\delta Q(t) -g(Q(t)) -\beta(Q(t),E(t))Q(t) \vspace{1ex}\\
&&+2e^{-\gamma\tau}\beta(Q(t-\tau),E(t-\tau))Q(t-\tau),\vspace{1ex}\\
\displaystyle\frac{dM}{dt}&=&-\mu M(t) +g(Q(t)),\vspace{1ex}\\
\displaystyle\frac{dE}{dt}&=&-k E(t)+f(M(t)).
\end{array}\right.
\end{equation}

First, one can notice that the solutions $Q(t)$, $M(t)$ and $E(t)$ of system (\ref{system}) are
nonnegative.

Let us suppose, by contradiction, that there exists $t_0>0$ and $\varepsilon>0$ such that $Q(t)>0$
for $t<t_0$, $Q(t_0)=0$ and $Q(t)<0$ for $t\in(t_0,t_0+\varepsilon)$. Then, since $g(0)=0$,
\begin{displaymath}
\frac{dQ}{dt}(t_0)=2e^{-\gamma\tau}\beta(Q(t_0-\tau),E(t_0-\tau))Q(t_0-\tau)>0.
\end{displaymath}
This yields a contradiction. Hence $Q(t)$ is nonnegative. Since the functions $g$ and $f$ are
nonnegative, we similarly obtain that $M(t)$ and $E(t)$ are nonnegative.

Secondly, one can notice that the solutions of (\ref{system}) are bounded when
$\delta+g^{\prime}(0)>0$.

We first concentrate on the solution $E(t)$. Using a classical variation of constant formula, we
obtain, for $t\geq0$,
\begin{displaymath}
E(t)=e^{-kt}E(0)+e^{-kt}\int_0^te^{ks}f(M(s))ds.
\end{displaymath}
Since the function $f$ is decreasing and bounded, we have
\begin{displaymath}
E(t)\leq e^{-kt}E(0)+\frac{f(0)}{k}(1-e^{-kt}).
\end{displaymath}
This yields that, if $E(0)\leq f(0)/k$, then $E(t)\leq f(0)/k$ and, if $f(0)/k<E(0)$, then
$E(t)\leq E(0)$. Consequently, $E(t)$ is bounded.

Now we focus on the solution $Q(t)$. If $Q$ is bounded then the mapping $t\mapsto g(Q(t))$ is
bounded so we will obtain that $M(t)$ is bounded using similar arguments than for the above case.

Let $C>0$ be a bound of $E$ and assume that $\lim_{Q\to\infty}\beta(Q,E)=0$, for all $E\geq0$, and
that $\delta+g^{\prime}(0)>0$. If $2e^{-\gamma\tau}\beta(0,C)\geq\delta+g^{\prime}(0)$, then, since
the mapping $Q\mapsto\beta(Q,E)$ is decreasing and tends to zero at infinity, there exists
$Q_0\geq0$ such that $2e^{-\gamma\tau}\beta(Q_0,C)=\delta+g^{\prime}(0)$ and
$2e^{-\gamma\tau}\beta(Q,C)<\delta+g^{\prime}(0)$ for $Q>Q_0$. If
$2e^{-\gamma\tau}\beta(0,C)>\delta+g^{\prime}(0)$ then this result holds for $Q_0=0$.

We then set
\begin{displaymath}
Q_1 := 2e^{-\gamma\tau}\frac{\beta(0,C)Q_0}{\delta+g^{\prime}(0)}.
\end{displaymath}
Let $Q\geq Q_1$ be fixed and let $0\leq y\leq Q$. If $y\leq Q_0$, then
\begin{displaymath}
2e^{-\gamma\tau}\beta(y,C)y\leq 2e^{-\gamma\tau}\beta(0,C)Q_0 =(\delta+g^{\prime}(0))Q_1\leq
(\delta+g^{\prime}(0))Q.
\end{displaymath}
On the other hand, if $y>Q_0$, then
\begin{displaymath}
2e^{-\gamma\tau}\beta(y,C)y<(\delta+g^{\prime}(0))y < (\delta+g^{\prime}(0))Q.
\end{displaymath}
Thus,
\begin{displaymath}
2e^{-\gamma\tau} \max_{0\leq y\leq Q}\beta(y,C)y \leq (\delta+g^{\prime}(0))Q, \quad \textrm{ for }
Q\geq Q_1.
\end{displaymath}
We assume now, by contradiction, that $\limsup Q(t)=+\infty$. Then there exists $t_0>\tau$ such
that
\begin{displaymath}
Q(t)\leq Q(t_0), \ \textrm{ for } t\in[t_0-\tau,t_0], \quad \textrm{ and } \quad Q(t_0)>Q_1.
\end{displaymath}
Since the function $E\mapsto\beta(Q,E)$ is increasing, we deduce, from (\ref{propg}), that
\begin{displaymath}
Q^{\prime}(t_0)\leq g^{\prime}(0)Q(t_0)-g(Q(t_0)) - \beta(Q(t_0),E(t_0))Q(t_0) <0.
\end{displaymath}
We obtain a contradiction so $Q$ is bounded.

These results are stated in the following proposition.

\begin{proposition}
Assume that $\lim_{Q\to\infty}\beta(Q,E)=0$, for all $E\geq0$, and that $\delta+g^{\prime}(0)>0$.
Then the solutions of system (\ref{system}) are bounded.
\end{proposition}

A solution $(\overline{Q},\overline{M},\overline{E})$ of (\ref{system}) is a steady-state or
equilibrium if
$$\frac{d\overline{Q}}{dt}=\frac{d\overline{M}}{dt}=\frac{d\overline{E}}{dt}=0,$$
that is
\begin{equation}\label{steadystate}
\left\{\begin{array}{rcl}
\delta\overline{Q}+g(\overline{Q})&=&(2e^{-\gamma\tau}-1)\beta(\overline{Q},\overline{E})\overline{Q},\vspace{1ex}\\
\mu\overline{M}&=&g(\overline{Q}),\vspace{1ex}\\
k\overline{E}&=&f(\overline{M}).
\end{array}\right.
\end{equation}
We make some remarks. It would be nonsense to suppose that the rates $\mu$ and $k$ may vanish,
because we cannot allow the blood cell population to grow indefinitely and growth factor is
necessarily degraded while in blood. Hence, we assume that $\mu>0$ and $k>0$.

Since $g(0)=0$, it follows that $(0,0,f(0)/k)$ is always a steady-state of (\ref{system}) that we
will call in the following the {\it trivial equilibrium} of (\ref{system}). This steady-state
corresponds to the extinction of the cell population with a saturation of the growth factor
concentration.

From now on and in the following, we assume that
\begin{equation}\label{hypbeta}
\lim_{Q\to+\infty}\beta\left(Q,\frac{1}{k}f\left(\frac{1}{\mu}g(Q)\right)\right)=0.
\end{equation}
Since $\lim_{Q\to+\infty}g(Q)=+\infty$ and $\lim_{M\to+\infty}f(M)=0$, this property holds true if
$\beta(Q,0)=0$, for all $Q\geq0$, or $\lim_{Q\to+\infty}\beta(Q,E)=0$, for all $E\geq0$.

Let us assume that (\ref{system}) has a nontrivial positive steady-state $(Q^*,M^*,E^*)$, that is,
$Q^*>0$, $M^*>0$ and $E^*>0$ satisfy (\ref{steadystate}). Then
\begin{equation}\label{steadystate2}
\left(2e^{-\gamma\tau}-1\right)\beta(Q^*,E^*)\!=\!\delta+\displaystyle\frac{g(Q^*)}{Q^*},\quad
M^*\!=\!\displaystyle\frac{1}{\mu}g(Q^*)\quad\! \textrm{and} \quad\!
E^*\!=\!\displaystyle\frac{1}{k}f(M^*).
\end{equation}
Necessarily, we have
\begin{displaymath}
2e^{-\gamma\tau}-1>0,
\end{displaymath}
which is equivalent to
\begin{equation}\label{tau}
\tau<\frac{\ln(2)}{\gamma}.
\end{equation}
We assume that (\ref{tau}) holds. Since $Q^*>0$ and
\begin{displaymath}
E^*=\frac{1}{k}f(M^*)=\frac{1}{k}f\left(\frac{1}{\mu}g(Q^*)\right),
\end{displaymath}
we must have
\begin{equation}\label{eqbeta}
\left(2e^{-\gamma\tau}-1\right)\beta\left(Q^*,\frac{1}{k}f\left(\frac{1}{\mu}g(Q^*)\right)\right)=\delta
+\displaystyle\frac{g(Q^*)}{Q^*} .
\end{equation}
Let us define
\begin{equation}\label{betatilda}
\widetilde{\beta}(Q) := \beta\left(Q,\frac{1}{k}f\left(\frac{1}{\mu}g(Q)\right)\right).
\end{equation}
Since $f$ is decreasing, $g$ is nondecreasing and $\beta(Q,E)$ is decreasing with respect to $Q$
and increasing with respect to $E$, we obtain
\begin{displaymath}
\widetilde{\beta}^{\prime}(Q)\!=\!\frac{\partial\beta}{\partial
Q}\left(Q,\frac{1}{k}f\left(\frac{1}{\mu}g(Q)\right)\right) +
\frac{g^{\prime}(Q)}{k\mu}f{^\prime}\left(\frac{1}{\mu}g(Q)\right)\frac{\partial\beta}{\partial
E}\left(Q,\frac{1}{k}f\left(\frac{1}{\mu}g(Q)\right)\right)<0.
\end{displaymath}
Therefore, $\widetilde{\beta}$ is decreasing,
\begin{displaymath}
\widetilde{\beta}(0)= \beta\left(0,\frac{1}{k}f(0)\right) \qquad \textrm{ and } \qquad
\lim_{Q\to+\infty}\widetilde{\beta}(Q)=0.
\end{displaymath}
Moreover, from (\ref{propg}), the function $Q\mapsto \delta+g(Q)/Q$ is nondecreasing. Consequently,
equation (\ref{eqbeta}) has a positive solution, which is unique, if and only if
\begin{equation}\label{condexist1}
\delta+g^{\prime}(0)<\left(2e^{-\gamma\tau}-1\right)\beta\left(0,\frac{1}{k}f(0)\right).
\end{equation}
These results are summarized in the next proposition.

\begin{proposition}\label{propequilibria}
Assume that $\mu>0$ and $k>0$.
\begin{description}
\item[(i)] If
\begin{equation}\label{condexist0}
\delta+g^{\prime}(0)>\left(2e^{-\gamma\tau}-1\right)\beta\left(0,\frac{1}{k}f(0)\right),
\end{equation}
then system (\ref{system}) has a unique steady-state $(0,0,f(0)/k);$
\item[(ii)] If condition (\ref{condexist1}) holds, then system (\ref{system}) has
two steady-states: a trivial one $(0,0,f(0)/k)$ and a nontrivial positive one $(Q^*,M^*,E^*)$,
where $Q^*$ is the unique positive solution of (\ref{eqbeta}), $M^*=g(Q^*)/\mu$ and $E^*=f(M^*)/k$.
\end{description}
\end{proposition}

The above proposition indicates that system (\ref{system}) undergoes a transcritical bifurcation
when $\delta+g^{\prime}(0)=\left(2e^{-\gamma\tau}-1\right)\beta(0,f(0)/k)$. Since the trivial
steady-state $(0,0,f(0)/k)$ corresponds, biologically, to the extinction of the cell population and
a saturation of the growth factor concentration, it is not the more interesting equilibrium. It
describes a pathological situation that can only lead to death without appropriate treatment.
Therefore we focus on the local stability analysis of the other steady-state $(Q^*,M^*,E^*)$. One
can check that condition (\ref{condexist1}), which ensures the existence of this steady-state, is
equivalent to
\begin{equation}\label{condexist2}
\delta\!+\!g^{\prime}(0)\!<\!\beta\left(0,\frac{f(0)}{k}\right)\!\quad\! \textrm{and} \quad\!
0\!\leq\!\tau\!<\!\tau_{max}\!:=\!\frac{1}{\gamma}\!\ln\!\left(\!\frac{2\beta\left(0,\frac{f(0)}{k}\right)}
{\delta\!+\!g^{\prime}(0)\!+\!\beta\left(0,\frac{f(0)}{k}\right)}\!\right).
\end{equation}

One can notice that, using the implicit function theorem, we can easily show that the steady-states
$Q^*$, $M^*$ and $E^*$ are continuously differentiable with respect to the cell cycle duration
$\tau$. Moreover, $Q^*$ and $M^*$ are decreasing functions of $\tau$ and $E^*$ is an increasing
function of $\tau$, with $\lim_{\tau\to \tau_{max}}(Q^*(\tau),M^*(\tau),E^*(\tau))=(0,0,f(0)/k)$.

\section{Local Stability Analysis} \label{sectionstability}

We concentrate on the study of the stability of the nontrivial equilibrium $(Q^*,M^*,E^*)$. Hence,
we assume throughout this section that $\mu, k>0$ and condition (\ref{condexist2}) holds.

We investigate the local asymptotic stability of the steady-state $(Q^*,M^*,E^*)$. The delay is
often seen as a destabilization parameter (see, for example, Metz and Diekmann \cite{md}, Fortin
and Mackey \cite{fm1999}), so we perform this stability analysis with respect to the delay
parameter $\tau$, which represents the cell cycle duration. We recall that $Q^*$, $M^*$ and $E^*$
satisfy (\ref{steadystate2}).

We linearize (\ref{system}) around the equilibrium $(Q^*,M^*,E^*)$. Set
\begin{displaymath}
q(t)=Q(t)-Q^*, \quad m(t)=M(t)-M^* \quad \textrm{ and } \quad e(t)=E(t)-E^*.
\end{displaymath}
The linearized system is
\begin{equation}\label{linsystem}
\left\{\begin{array}{rcl}
\displaystyle\frac{dq}{dt}&=&-Aq(t)+Bq(t-\tau)-Ce(t)+De(t-\tau),\vspace{1ex}\\
\displaystyle\frac{dm}{dt}&=&-\mu m(t) +G q(t),\vspace{1ex}\\
\displaystyle\frac{de}{dt}&=&-k e(t)-Hm(t),
\end{array}\right.
\end{equation}
where the real coefficients $A, B, C, D, G$ and $H$ are defined by
\begin{equation}\label{coefAB}
\begin{array}{rcl}
A&=&\delta+g^{\prime}(Q^*)+\beta(Q^*,E^*)+\beta^{\prime}_1(Q^*,E^*)Q^*,\vspace{1ex}\\
B&=&2e^{-\gamma\tau}\left[\beta(Q^*,E^*)+\beta^{\prime}_1(Q^*,E^*)Q^*\right],\vspace{1ex}\\
C&=&\beta^{\prime}_2(Q^*,E^*)Q^*>0,\vspace{1ex}\\
D&=&2e^{-\gamma\tau}\beta^{\prime}_2(Q^*,E^*)Q^*>0,\vspace{1ex}\\
G&=&g^{\prime}(Q^*)>0,\vspace{1ex}\\
H&=&-f^{\prime}(M^*)>0.
\end{array}
\end{equation}
One can notice that these coefficients depend, explicitly or implicitly, on the parameter $\tau$
through the equilibrium values $Q^*$, $M^*$ and $E^*$. However, we do not stress this dependence
when we write the coefficients. Moreover, from the assumptions on $\beta$, $g$ and $f$, the
coefficients $C$, $D$, $G$ and $H$ are strictly positive.

In the above definitions, we have used the notations
\begin{displaymath}
\beta^{\prime}_1(Q,E):=\frac{\partial\beta}{\partial Q}(Q,E) \qquad \textrm{ and } \qquad
\beta^{\prime}_2(Q,E):=\frac{\partial\beta}{\partial E}(Q,E).
\end{displaymath}
Furthermore, one can notice that
\begin{displaymath}
A-B=g^{\prime}(Q^*)-\frac{g(Q^*)}{Q^*}-\left(2e^{-\gamma\tau}-1\right)\beta^{\prime}_1(Q^*,E^*)Q^*\geq0
\end{displaymath}
and
\begin{displaymath}
D-C=\left(2e^{-\gamma\tau}-1\right)\beta^{\prime}_2(Q^*,E^*)Q^*>0.
\end{displaymath}

System (\ref{linsystem}) can be written in the matrix form
\begin{displaymath}
\frac{dX}{dt}=\mathcal{A}_1X(t)+\mathcal{A}_2X(t-\tau)
\end{displaymath}
with
\begin{displaymath}
X(t)=\left(\begin{array}{c}
q(t)\vspace{1ex}\\
m(t)\vspace{1ex}\\
e(t)
\end{array}\right), \quad
\mathcal{A}_1=\left(\begin{array}{ccc}
-A&0&-C\vspace{1ex}\\
G&-\mu&0\vspace{1ex}\\
0&-H&-k\end{array}\right) \quad \textrm{ and } \quad \mathcal{A}_2=\left(\begin{array}{ccc}
B&0&D\vspace{1ex}\\
0&0&0\vspace{1ex}\\
0&0&0\end{array}\right).
\end{displaymath}
Consequently, the characteristic equation of (\ref{linsystem}) is given by
\begin{displaymath}
\det\left(\lambda I-\mathcal{A}_1-\mathcal{A}_2e^{-\lambda\tau} \right)=0,
\end{displaymath}
which reduces to
\begin{equation}\label{ce}
(\lambda+\mu)(\lambda+k)\left(\lambda+A-Be^{-\lambda\tau}\right)-G
H\left(C-De^{-\lambda\tau}\right)=0.
\end{equation}
We recall the following result: the trivial solution of system (\ref{linsystem}), or equivalently
the steady-state of system (\ref{system}), is asymptotically stable if all roots of (\ref{ce}) have
negative real parts, and the stability is lost only if characteristic roots cross the axis from
left to right, or right to left, that is if purely imaginary roots appear.

\begin{remark}{\rm
If we linearize system (\ref{system}) around its trivial steady-state $(0,0,f(0)/k)$, we obtain
\begin{displaymath}
\begin{array}{rclcrcl}
A&=&\delta+g^{\prime}(0)+\beta(0,f(0)/k)>0,&\ &D&=&0,\vspace{1ex}\\
B&=&2e^{-\gamma\tau}\beta(0,f(0)/k)>0,&\ &G&=&g^{\prime}(0)>0,\vspace{1ex}\\
C&=&0,\vspace{1ex}&\ &H&=&-f^{\prime}(0)>0.
\end{array}
\end{displaymath}
Therefore, the characteristic equation (\ref{ce}) of the linearized system, around the trivial
steady-state, becomes
\begin{displaymath}
(\lambda+\mu)(\lambda+k)\left(\lambda+A-Be^{-\lambda\tau}\right)=0.
\end{displaymath}
It follows that $\lambda=-\mu<0$, $\lambda=-k<0$ or
\begin{equation}\label{ce4}
\lambda+A-Be^{-\lambda\tau}=0.
\end{equation}
Studying the sign of the real parts of the roots of (\ref{ce4}), we obtain the following
proposition which deals with the local asymptotic stability of the trivial steady-state of
(\ref{system}).}
\end{remark}

\begin{proposition}
Assume that $\mu>0$ and $k>0$. If condition (\ref{condexist0}) holds, then the trivial steady-state
of system (\ref{system}) is locally asymptotically stable for all $\tau\geq0$, and if condition
(\ref{condexist1}) holds, then it is unstable for all $\tau\geq0$.
\end{proposition}

\begin{proof}
First notice that, when $\tau=0$, $\lambda=B-A$ so $\lambda>0$ if condition (\ref{condexist1})
holds and $\lambda<0$ if condition (\ref{condexist0}) holds.

Let $\tau>0$ be fixed. Setting $\nu=\lambda\tau$, the characteristic equation (\ref{ce4}) is
equivalent to
\begin{displaymath}
(\nu+A\tau)e^{\nu}-B\tau=0.
\end{displaymath}
From \cite{hayes}, we know that $\textrm{Re}(\nu)<0$ if and only if
\begin{displaymath}
A\tau>-1, \quad A\tau-B\tau>0, \quad \textrm{ and } \quad B\tau<\zeta\sin(\zeta)-A\tau\cos(\zeta),
\end{displaymath}
where $\zeta$ is the unique solution of
\begin{displaymath}
\zeta=-A\tau\tan(\zeta), \qquad \zeta\in(0,\pi).
\end{displaymath}
Since $A>0$ and $\tau>0$, condition $A\tau>-1$ is satisfied.

If we assume that condition (\ref{condexist0}) holds, then $A>B$ so $A\tau-B\tau>0$. By
contradiction, suppose that $B\tau>\zeta\sin(\zeta)-A\tau\cos(\zeta)$. Then, from the definition of
$\zeta$,
\begin{displaymath}
B\tau>-\frac{A\tau}{\cos(\zeta)}.
\end{displaymath}
Since $A>B>0$, it follows that
\begin{displaymath}
1>-\frac{1}{\cos(\zeta)}.
\end{displaymath}
Consequently, $\cos(\zeta)>0$ and $\zeta\in(\pi/2,\pi)$. We deduce that $\tan(\zeta)>0$ so
\begin{displaymath}
-A\tau\tan(\zeta)<0<\zeta.
\end{displaymath}
This gives a contradiction. Therefore $B\tau<\zeta\sin(\zeta)-A\tau\cos(\zeta)$, and all roots of
(\ref{ce4}) have negative real parts. The trivial steady-state is then locally asymptotically
stable for all $\tau>0$.

Assume now that condition (\ref{condexist1}) holds. Then $A<B$ and $A\tau-B\tau<0$. Consequently,
for all $\tau>0$, (\ref{ce4}) has roots with nonnegative real parts and the trivial steady-state is
unstable.
\end{proof}

The above results indicate that the trivial steady-state of (\ref{system}) is locally
asymptotically stable when it is the only equilibrium and unstable as soon as the nontrivial
equilibrium exists. When the transcritical bifurcation occurs, that is when the two steady-states
coincide, the trivial steady-state is stable.

We now return to the analysis of the local asymptotic stability of the nontrivial steady-state
$(Q^*,M^*,E^*)$ of system (\ref{system}).

Equation (\ref{ce}) takes the general form
\begin{equation}\label{ce1}
P(\lambda,\tau)+Q(\lambda,\tau)e^{-\lambda\tau}=0
\end{equation}
with
\begin{eqnarray}
P(\lambda,\tau)&=&\lambda^3+a_1(\tau)\lambda^2+a_2(\tau)\lambda+a_3(\tau),\nonumber\\
Q(\lambda,\tau)&=&a_4(\tau)\lambda^2+a_5(\tau)\lambda+a_6(\tau),\nonumber
\end{eqnarray}
where
\begin{displaymath}
\begin{array}{rclcrcl}
a_1(\tau)&=&\mu+k+A,&&a_4(\tau)&=&-B,\vspace{1ex}\\
a_2(\tau)&=&\mu k +A(\mu+k),&&a_5(\tau)&=&-B(\mu+k),\vspace{1ex}\\
a_3(\tau)&=&\mu kA-G HC,&&a_6(\tau)&=&-B\mu k+G HD.
\end{array}
\end{displaymath}
We can check that, for all $\tau\in[0,\tau_{max})$,
\begin{displaymath}
a_1(\tau)+a_4(\tau)=\mu+k+A-B>0,
\end{displaymath}
\begin{displaymath}
a_2(\tau)+a_5(\tau)=\mu k +(A-B)(\mu+k)>0,
\end{displaymath}
and
\begin{displaymath}
a_3(\tau)+a_6(\tau)=\mu k(A-B)+G H(D-C)>0.
\end{displaymath}
We will remember, in the following, that
\begin{equation}\label{propa}
a_i(\tau)+a_{i+3}(\tau)>0 \qquad \textrm{ for } i=1,2,3.
\end{equation}

Let us examine the case $\tau=0$. This case is of importance, because it can be necessary that the
nontrivial positive steady-state of (\ref{system}) is stable when $\tau=0$ to be able to obtain the
local stability for all nonnegative values of the delay, or to find a critical value which could
destabilize the steady-state (see Theorem \ref{theorembif}).

When $\tau=0$, the characteristic equation (\ref{ce1}) reduces to
\begin{equation}\label{ce2}
\lambda^3+[a_1(0)+a_4(0)]\lambda^2+[a_2(0)+a_5(0)]\lambda+a_3(0)+a_6(0)=0.
\end{equation}
The Routh-Hurwitz criterion says that all roots of (\ref{ce2}) have negative real parts if and only
if
\begin{displaymath}
\begin{array}{l}
a_1(0)+a_4(0)>0,\vspace{1ex}\\
a_3(0)+a_6(0)>0,
\end{array}
\end{displaymath}
and
\begin{equation}\label{condstabzero}
[a_1(0)+a_4(0)][a_2(0)+a_5(0)]>a_3(0)+a_6(0).
\end{equation}
From (\ref{propa}), it follows that all characteristic roots of (\ref{ce2}) have negative real
parts if and only if (\ref{condstabzero}) holds.

\begin{proposition}\label{propstabtauzero}
When $\tau=0$, the nontrivial steady-state $(Q^*,M^*,E^*)$ of (\ref{system}) is locally
asymptotically stable if and only if
\begin{equation}\label{stabzero}
(\mu+k)\left[\mu k+(A-B)(\mu+k+A-B)\right]>G H(D-C),
\end{equation}
where $A$, $B$, $C$, $D$, $G$ and $H$ are given by (\ref{coefAB}).
\end{proposition}

In the following, we investigate the existence of purely imaginary roots $\lambda=i\omega$,
$\omega\in\mathbb{R}$, of (\ref{ce1}). Equation (\ref{ce1}) takes the form of a third degree
exponential polynomial in $\lambda$. In 2001, Ruan and Wey \cite{rw2001} gave sufficient conditions
for the existence of zeros for such an equation, but only in the case where the coefficients of the
polynomial functions $P$ and $Q$ do not depend on the delay $\tau$, that is when the characteristic
equation (\ref{ce1}) is given by $P(\lambda)+Q(\lambda)e^{-\lambda\tau}=0$. Since all the
coefficients of $P$ and $Q$ depend on $\tau$, we cannot apply their results directly. In 2002,
however, Beretta and kuang \cite{bk2002} established a geometrical criterion which gives the
existence of purely imaginary roots for a characteristic equation with delay dependent
coefficients. We are going to apply this criterion to equation (\ref{ce1}) in order to obtain
stability results for equation (\ref{linsystem}). In the following, we use the same notations as in
Beretta and Kuang \cite{bk2002} to make things easier for the reader.

We first have to verify the following properties, for all $\tau\in[0,\tau_{max})$:
\begin{description}
\item[(i)] $P(0,\tau)+Q(0,\tau)\neq0$;
\item[(ii)] $P(i\omega,\tau)+Q(i\omega,\tau)\neq0$;
\item[(iii)] $\limsup\left\{\left|\frac{Q(\lambda,\tau)}{P(\lambda,\tau)}\right|;
|\lambda|\to\infty, \textrm{Re}\lambda\geq0\right\}<1$;
\item[(iv)] $F(\omega,\tau) := |P(i\omega,\tau)|^2-|Q(i\omega,\tau)|^2$ has a finite number of
zeros.
\end{description}
Properties (i), (ii) and (iii) can be easily verified. Let $\tau\in[0,\tau_{max})$. Using
(\ref{propa}), a simple computation gives
\begin{displaymath}
P(0,\tau)+Q(0,\tau)=a_3(\tau)+a_6(\tau)>0.
\end{displaymath}
Moreover,
\begin{displaymath}
\begin{array}{rcl}
P(i\omega,\tau)+Q(i\omega,\tau)&=&\left[-(a_1(\tau)+a_4(\tau))\omega^2+a_3(\tau)+a_6(\tau)\right]\\
&&+i\left[-\omega^3+(a_2(\tau)+a_5(\tau))\omega\right],
\end{array}
\end{displaymath} so (ii) is true.
Finally,
\begin{displaymath}
\left|\frac{Q(\lambda,\tau)}{P(\lambda,\tau)}\right|\underset{|\lambda|\to\infty}{\sim}
\left|\frac{a_4(\tau)}{\lambda}\right|,
\end{displaymath}
therefore (iii) is also true.

Now, let $F$ be defined as in (iv). Since
\begin{displaymath}
|P(i\omega,\tau)|^2=\omega^6+\left[a_1^2(\tau)-2a_2(\tau)\right]\omega^4+\left[a_2^2(\tau)-2a_1(\tau)a_3(\tau)\right]\omega^2+a_3^2(\tau)
\end{displaymath}
and
\begin{displaymath}
|Q(i\omega,\tau)|^2=a_4^2(\tau)\omega^4+\left[a_5^2(\tau)-2a_4(\tau)a_6(\tau)\right]\omega^2+a_6^2(\tau),
\end{displaymath}
we have
\begin{displaymath}
F(\omega,\tau)=\omega^6+b_1(\tau)\omega^4+b_2(\tau)\omega^2+b_3(\tau)
\end{displaymath}
with
\begin{displaymath}
\begin{array}{rcl}
b_1(\tau)&=&a_1^2(\tau)-2a_2(\tau)-a_4^2(\tau),\vspace{1ex}\\
b_2(\tau)&=&a_2^2(\tau)+2a_4(\tau)a_6(\tau)-2a_1(\tau)a_3(\tau)-a_5^2(\tau),\vspace{1ex}\\
b_3(\tau)&=&a_3^2(\tau)-a_6^2(\tau).
\end{array}
\end{displaymath}
One can check that
\begin{displaymath}
b_1(\tau)=\mu^2+k^2+A^2-B^2,
\end{displaymath}
and
\begin{displaymath}
\begin{array}{rcl}
b_2(\tau)&=&\mu^2k^2+(A^2-B^2)(\mu^2+k^2)+2G H\left[C(\mu+k+A)-BD\right],\vspace{1ex}\\
b_3(\tau)&=&\mu^2k^2(A^2-B^2)+G^2H^2(C^2-D^2)+2\mu kGH(BD-AC),
\end{array}
\end{displaymath}
where $A$, $B$, $C$, $D$, $G$ and $H$ are given by (\ref{coefAB}). It is obvious that property (iv)
is satisfied.

Now assume that $\lambda=i\omega$, $\omega\in\mathbb{R}$, is a purely imaginary characteristic root
of (\ref{ce1}). Separating real and imaginary parts, we can show that $(\omega,\tau)$ satisfies
\begin{eqnarray}
-a_1(\tau)\omega^2+a_3(\tau)&=&-\left[-a_4(\tau)\omega^2+a_6(\tau)\right]\cos(\omega\tau)-a_5(\tau)\omega\sin(\omega\tau),\vspace{1ex}\label{e1}\\
-\omega^3+a_2(\tau)\omega&=&-a_5(\tau)\omega\cos(\omega\tau)+\left[-a_4(\tau)\omega^2+a_6(\tau)\right]\sin(\omega\tau)\label{e2}.
\end{eqnarray}
One can check that, if $(\omega,\tau)$ is a solution of system (\ref{e1})--(\ref{e2}), then
$(-\omega,\tau)$ is also a solution of (\ref{e1})--(\ref{e2}). Hence, if $i\omega$ is a purely
imaginary characteristic root of (\ref{ce1}), its conjugate has the same property. Consequently, we
only look in the following for purely imaginary roots of (\ref{ce1}) with positive imaginary part.

System (\ref{e1})--(\ref{e2}) yields
\begin{eqnarray}
\cos(\omega\tau)&=&\displaystyle\frac{\left(a_5-a_1a_4\right)\omega^4+\left(a_1a_6+a_3a_4-a_2a_5\right)\omega^2-a_3a_6}{a_4^2\omega^4+\left(a_5^2-2a_4a_6\right)\omega^2+a_6^2},\label{e3}\\
\sin(\omega\tau)&=&\displaystyle\frac{a_4\omega^5+\left(a_1a_5-a_2a_4-a_6\right)\omega^3+\left(a_2a_6-a_3a_5\right)\omega}{a_4^2\omega^4+\left(a_5^2-2a_4a_6\right)\omega^2+a_6^2}\label{e4},
\end{eqnarray}
where we deliberately omit the dependence of the $a_i$ on $\tau$.

A necessary condition for this system to have solutions is that the sum of the squares of the right
hand side terms equals one. By remarking that system (\ref{e3})--(\ref{e4}) can be written
\begin{displaymath}
\cos(\omega\tau)=\textrm{Im}\left(\frac{P(i\omega,\tau)}{Q(i\omega,\tau)}\right) \quad \textrm{ and
} \quad \sin(\omega\tau)=-\textrm{Re}\left(\frac{P(i\omega,\tau)}{Q(i\omega,\tau)}\right),
\end{displaymath}
then this condition is
\begin{displaymath}
|P(i\omega,\tau)|^2=|Q(i\omega,\tau)|^2.
\end{displaymath}
That is
\begin{displaymath}
F(\omega,\tau)=0.
\end{displaymath}
The polynomial function $F$ can be written as
\begin{displaymath}
F(\omega,\tau)=h(\omega^2,\tau),
\end{displaymath}
where $h$ is a third degree polynomial, defined by
\begin{equation}\label{h}
h(z,\tau):=z^3+b_1(\tau)z^2+b_2(\tau)z+b_3(\tau).
\end{equation}
We set
\begin{equation}\label{Delta}
\Delta(\tau)=b_1^2(\tau)-3b_2(\tau),
\end{equation}
and, when $\Delta(\tau)\geq0$,
\begin{equation}\label{z0}
z_0(\tau)=\frac{-b_1(\tau)+\sqrt{\Delta(\tau)}}{3}.
\end{equation}
We then have the following lemma (see \cite{rw2001}).

\begin{lemma}\label{lemerp}
Let $\tau\in[0,\tau_{max})$ and $\Delta(\tau)$ and $z_0(\tau)$ be defined by (\ref{Delta}) and
(\ref{z0}), respectively. Then $h(\cdot,\tau)$, defined in (\ref{h}), has positive roots if and
only if
\begin{equation}\label{cond2}
b_3(\tau)<0 \quad \textrm{ or } \quad b_3(\tau)\geq0, \ \Delta(\tau)\geq0, \ z_0(\tau)>0 \
\textrm{and} \ h(z_0(\tau),\tau)<0.
\end{equation}
\end{lemma}

\begin{proof}
Details of the proof are given in \cite{rw2001}, Lemma 2.1.
\end{proof}

Conditions $\Delta(\tau)\geq0$, $z_0(\tau)>0$ and $h(z_0(\tau),\tau)<0$ cannot be easily checked.
In the following lemma we express them using the coefficients $b_i$, $i=1,2,3$.

\begin{lemma}\label{b3pos}
Let $\tau\geq0$ be such that $b_3(\tau)\geq0$. Then $\Delta(\tau)\geq0$, $z_0(\tau)>0$ and
$h(z_0(\tau),\tau)<0$ if and only if
\begin{displaymath}
\begin{array}{cl}
(i)& b_2(\tau)<0 \quad \textrm{ or } \quad b_1(\tau)<0\leq
b_2(\tau)<\displaystyle\frac{b_1^2(\tau)}{3}, \textrm{ and} \vspace{1ex}\\
(ii)& 2\Delta(\tau)z_0(\tau)+b_1(\tau)b_2(\tau)-9b_3(\tau)>0.
\end{array}
\end{displaymath}
\end{lemma}

\begin{proof}
Let $\tau$ be given such that $b_3(\tau)\geq0$. We do not mention, in the following, the dependence
of the coefficients $b_i$ on $\tau$.

We have
\begin{displaymath}
\Delta \geq0 \qquad \textrm{ if and only if } \qquad b_1^2\geq 3b_2.
\end{displaymath}
If $b_2<0$, this result holds true. Otherwise, it is necessary that $b_1^2\geq 3b_2$. In this
latter case, if $b_1<0$, then $z_0>0$ and, if $b_1\geq0$, then $z_0>0$ if and only if $b_2<0$.
Therefore $z_0>0$ if and only if
\begin{equation}\label{nc1}
b_2<0 \qquad \textrm{ or } \qquad b_1<0\leq b_2<\frac{b_1^2}{3}.
\end{equation}
Under the assumption (\ref{nc1}), $h^{\prime}$, given by
\begin{displaymath}
h^{\prime}(z)=3z^2+2b_1z+b_2,
\end{displaymath}
has two roots,
\begin{displaymath}
z_-=-\frac{1}{3}(b_1+d) \quad \textrm{ and } \quad z_+=-\frac{1}{3}(b_1-d)
\end{displaymath}
with $z_-<z_+$ and $d=\sqrt{b_1^2-3b_2}$ (in fact $z_+=z_0>0$). A simple computation gives
\begin{displaymath}
h(z_+)=\frac{2}{27}(b_1^3-d^3)-\frac{b_1b_2}{3}+b_3.
\end{displaymath}
Noticing that
\begin{displaymath}
b_1^3-d^3=(b_1-d)(2b_1^2-3b_2+b_1d)=-3z_+(b_1^2+b_1d+d^2),
\end{displaymath}
we obtain
\begin{displaymath}
h(z_+)<0 \quad \Leftrightarrow \quad \frac{2}{3}z_+(b_1^2+b_1d+d^2)+b_1b_2-3b_3>0.
\end{displaymath}
Moreover,
\begin{displaymath}
b_1^2+b_1d+d^2=d^2+b_1(b_1+d)=d^2-3b_1z_-.
\end{displaymath}
So
\begin{displaymath}
h(z_+)<0 \quad \Leftrightarrow \quad \frac{2}{3}d^2z_+ -2b_1z_+z_-+b_1b_2-3b_3>0.
\end{displaymath}
Since $z_+z_-=b_2/3$, we eventually obtain
\begin{displaymath}
h(z_+)<0 \quad \Leftrightarrow \quad 2d^2z_++b_1b_2-9b_3>0.
\end{displaymath}
This ends the proof.
\end{proof}

From the previous lemma, condition (\ref{cond2}) is equivalent to
\begin{equation}\label{cond3}
b_3(\tau)<0 \qquad \textrm{ or } \qquad b_3(\tau)\geq0 \textrm{ and }(i)\textrm{-}(ii) \textrm{
hold true.}
\end{equation}
Let us show on an example that condition (\ref{cond3}) is satisfied.

Note that $b_3$ can be expressed as
\begin{displaymath}
\begin{array}{rcl}
b_3(\tau)&=&\mu^2k^2(A-B)(A+B)+G^2H^2(C-D)(C+D)\\
&&+2\mu kGH(B(D-C)+C(B-A)),
\end{array}
\end{displaymath}
where $A$, $B$, $C$, $D$, $G$ and $H$ are defined by (\ref{coefAB}). Since $C-D<0$ and $B-A\leq0$,
then $b_3(\tau)<0$ if $A+B\leq0$ and $B\leq0$. Moreover, from the definition of $B$, it follows
that $B\leq0$ if $A+B\leq0$. Consequently, a sufficient condition for $b_3(\tau)<0$ is $A+B\leq0$.

Let us assume that $g$ is a linear function, given by $g(Q)=GQ$, with $G>0$, and that
$\beta(Q,E)=\beta_1(Q)\beta_2(E)$, with $\beta_1(Q)=1/(1+Q^n)$, $n>0$, and $\beta_2$ an increasing
function satisfying $\beta_2(0)=0$. Then, for $\tau=0$,
\begin{displaymath}
A+B=[4\beta_1(Q^*)+3\beta^{\prime}_1(Q^*)Q^*]\beta_2(E^*).
\end{displaymath}
Since $E^*>0$, $\beta_2(E^*)>0$ and $A+B\leq0$ if and only if
\begin{displaymath}
4\beta_1(Q^*)+3\beta^{\prime}_1(Q^*)Q^*=\frac{4+(4-3n)(Q^*)^n}{(1+(Q^*)^n)^2}\leq0,
\end{displaymath}
that is
\begin{displaymath}
n>\frac{4}{3} \qquad \textrm{ and } \qquad Q^*\geq\left(\frac{4}{3n-4}\right)^{1/n}.
\end{displaymath}

Let $\delta$ and $G$ be such that
\begin{displaymath}
\delta+G<\widetilde{\beta}(1)=\frac{1}{2}\beta_2\left(\frac{1}{k}f\left(\frac{G}{\mu}\right)\right),
\end{displaymath}
where $\widetilde{\beta}$ is defined by (\ref{betatilda}), and let $\overline{n}>4/3$ be the unique
solution of
\begin{displaymath}
\frac{3\overline{n}}{3\overline{n}-4}+\ln\left(\frac{4}{3\overline{n}-4}\right)=0.
\end{displaymath}
Then, for $n>\overline{n}\approx 6.12$,
\begin{displaymath}
\widetilde{\beta}\left(\left(\frac{4}{3n-4}\right)^{1/n}\right)>\widetilde{\beta}(1)>\delta+G,
\end{displaymath}
so $Q^*\geq(4/(3n-4))^{1/n}$ and it follows that $A+B\leq0$.

Consequently, $b_3(0)<0$ and, using the continuity of $b_3$ with respect to $\tau$, we deduce that
there exists $\overline{\tau}>0$ such that (\ref{cond3}) is verified for
$\tau\in[0,\overline{\tau})$.

When the reintroduction rate $\beta$ only depends on the growth factor concentration $E$, condition
(\ref{cond3}) is also satisfied for $\tau$ close to zero. This is numerically obtained in Section
\ref{snilpo}.

We set $I := [0,\overline{\tau})$ an interval in which (\ref{cond3}) is satisfied, with
$0<\overline{\tau}\leq\tau_{max}$. From the above remarks, we can find functions $\beta$, $g$ and
$f$, and parameters values such that $\overline{\tau}$ exists. For $\tau\in I$ there exists at
least $\omega=\omega(\tau)>0$ such that $F(\omega(\tau),\tau)=0$.

Then, let $\theta(\tau)\in[0,2\pi]$ be defined for $\tau\in I$ by
\begin{eqnarray}
\cos(\theta(\tau))&=&\displaystyle\frac{\left(a_5-a_1a_4\right)\omega^4+\left(a_1a_6+a_3a_4-a_2a_5\right)\omega^2-a_3a_6}{a_4^2\omega^4+\left(a_5^2-2a_4a_6\right)\omega^2+a_6^2},\nonumber\\
\sin(\theta(\tau))&=&\displaystyle\frac{a_4\omega^5+\left(a_1a_5-a_2a_4-a_6\right)\omega^3+\left(a_2a_6-a_3a_5\right)\omega}{a_4^2\omega^4+\left(a_5^2-2a_4a_6\right)\omega^2+a_6^2},\nonumber
\end{eqnarray}
where $\omega=\omega(\tau)$. Since $F(\omega(\tau),\tau)=0$ for $\tau\in I$, it follows that
$\theta$ is well and uniquely defined for all $\tau\in I$.

One can check that $i\omega^*$, with $\omega^*=\omega(\tau^*)>0$ is a purely imaginary
characteristic root of (\ref{ce1}) if and only if $\tau^*$ is a root of the function $S_n$, defined
by
\begin{displaymath}
S_n(\tau)=\tau-\frac{\theta(\tau)+2n\pi}{\omega(\tau)}, \quad \tau\in I, \qquad \textrm{ with }
n\in\mathbb{N}.
\end{displaymath}
The following theorem is due to Beretta and Kuang \cite{bk2002}.

\begin{theorem}\label{theorembk}
Assume that the function $S_n(\tau)$ has a positive root $\tau^*\in I$ for some $n\in\mathbb{N}$.
Then a pair of simple purely imaginary roots $\pm i\omega(\tau^*)$ of (\ref{ce1}) exists at
$\tau=\tau^*$ and
\begin{equation}\label{tc1}
\textrm{{\em sign}}\left\{\frac{dRe(\lambda)}{d\tau}\bigg|_{\lambda=i\omega(\tau^*)}\right\}=
\textrm{{\em sign}}\left\{\frac{\partial F}{\partial \omega}(\omega(\tau^*),\tau^*)\right\}
\textrm{{\em sign}}\left\{\frac{dS_n(\tau)}{d\tau}\bigg|_{\tau=\tau^*}\right\}.
\end{equation}
\end{theorem}

\vspace{2ex}

Since
\begin{displaymath}
\frac{\partial F}{\partial \omega}(\omega,\tau)=2\omega\frac{\partial h}{\partial
z}(\omega^2,\tau),
\end{displaymath}
condition (\ref{tc1}) is equivalent to
\begin{displaymath}
\textrm{ sign}\left\{\frac{dRe(\lambda)}{d\tau}\bigg|_{\lambda=i\omega(\tau^*)}\right\}=
\textrm{sign}\left\{\frac{\partial h}{\partial z}(\omega^2(\tau^*),\tau^*)\right\} \textrm{
sign}\left\{\frac{dS_n(\tau)}{d\tau}\bigg|_{\tau=\tau^*}\right\}.
\end{displaymath}

We can easily observe that $S_n(0)<0$. Moreover, for all $\tau\in I$, $S_n(\tau)>S_{n+1}(\tau)$,
with $n\in\mathbb{N}$. Therefore, if $S_0$ has no root in $I$, then the $S_n$ functions have no
root in $I$ and, if the function $S_n(\tau)$ has positive roots $\tau\in I$ for some
$n\in\mathbb{N}$, there exists at least one root satisfying
\begin{displaymath}
\frac{dS_n}{d\tau}(\tau)>0.
\end{displaymath}
Using Proposition \ref{propstabtauzero}, we can conclude the existence of a Hopf bifurcation as
stated in the next theorem.

\begin{theorem}\label{theorembif}
Assume that $\mu$, $k>0$, condition (\ref{condexist1}) is satisfied and (\ref{stabzero}) holds
true.
\begin{description}
\item[(i)] If the function $S_0(\tau)$ has no positive root in
$I$, then the steady-state $(Q^*,M^*,E^*)$ is locally asymptotically stable for all $\tau\geq0$.
\item[(ii)] If the function $S_0(\tau)$ has at least one positive root in $I,$ then there exists
$\tau^*\in I$ such that the steady-state $(Q^*,M^*,E^*)$ is locally asymptotically stable for
$0\leq\tau<\tau^*$ and becomes unstable for $\tau\geq\tau^*$, with a Hopf bifurcation occurring
when $\tau=\tau^*$, if and only if
\begin{displaymath}
\frac{\partial h}{\partial z}(\omega^2(\tau^*),\tau^*)>0.
\end{displaymath}
\end{description}
\end{theorem}

In the next section, we illustrate the results established in Theorem \ref{theorembif}. We show, in
particularly, that our model can exhibit long periods oscillations, compared to the delay, that can
be related to experimental observations in patients with periodic hematological diseases.

\section{Numerical Illustrations: Long Periods Oscillations}\label{snilpo}

Let us assume that the introduction of resting cells in the proliferating phase is only triggered
by the growth factor concentration $E(t)$, that is $\beta=\beta(E(t))$. This assumption is based on
hypothesis made by Belair {\it et al.} \cite{bmm1995} for an erythropoiesis model. It describes,
for example, the fact that the cell population may only react to external stimuli and cannot be
directly sensitive to its own size. We assume that $\beta$ is given by
\begin{displaymath}
\beta(E)=\beta_0\frac{E}{1+E}, \qquad \beta_0>0.
\end{displaymath}
The functions $g$ and $f$ are defined by
\begin{displaymath}
g(Q)=GQ \qquad \textrm{ with } G>0,
\end{displaymath}
and
\begin{displaymath}
f(M)=\frac{a}{1+KM^r}, \qquad a,K>0, r>1.
\end{displaymath}
This latter function often occurs in enzyme kinetics. It has been used by Mackey \cite{m1978,m1979}
to describe the rate of introduction in the proliferating phase and by Belair {\it et al.}
\cite{bmm1995} to define the feedback from the blood to the growth factor production.

With these choices for the functions $\beta$, $g$ and $f$, our model involves 10 parameters,
including the delay $\tau$. Most of the values of these parameters can be found in the literature.
This is the case for the stem cells mortality rates in the bone marrow, $\delta$ and $\gamma$,
which are given by Mackey \cite{m1978} or Pujo-Menjouet and Mackey \cite{pm2004}, and for the
coefficients of the function $f$ and the disappearance rate $k$, given by Belair {\it et al.}
\cite{bmm1995} and Mahaffy {\it et al.} \cite{mbm1998}. However, it is not so easy to determine
values for the other parameters.

Experimental data indicate similar values for $\mu$ and $G$, since blood cells are almost produced
at the same rate than they are lost. Then, we will choose
\begin{equation}\label{val1}
\mu=0.02 \textrm{ d}^{-1} \qquad \textrm{ and } \qquad G=0.04 \textrm{ d}^{-1}.
\end{equation}
The coefficient $\beta_0$ represents the maximum rate of introduction in the proliferating phase
and also the value of $\beta^{\prime}(0)$. It strongly depends on the nature of the growth factor.
Using data about erythropoiesis, we choose
\begin{equation}\label{val2}
\beta_0=0.5,
\end{equation}
which is less than the maximal rate of introduction proposed by Mackey \cite{m1978}, but seems
sufficiently large for erythropoiesis modelling.

Other parameters are given by (see \cite{bmm1995,m1978})
\begin{equation}\label{val3}
\begin{array}{c}
\delta=0.01 \textrm{ d}^{-1},  \ \gamma=0.2 \textrm{ d}^{-1}, \ k=2.8 \textrm{ d}^{-1},\\
a=6570, \quad K=0.0382 \ \textrm{ and } \ r=7.
\end{array}
\end{equation}

\begin{figure}[hpt]
\begin{center}
\includegraphics[width=9cm, height=7cm]{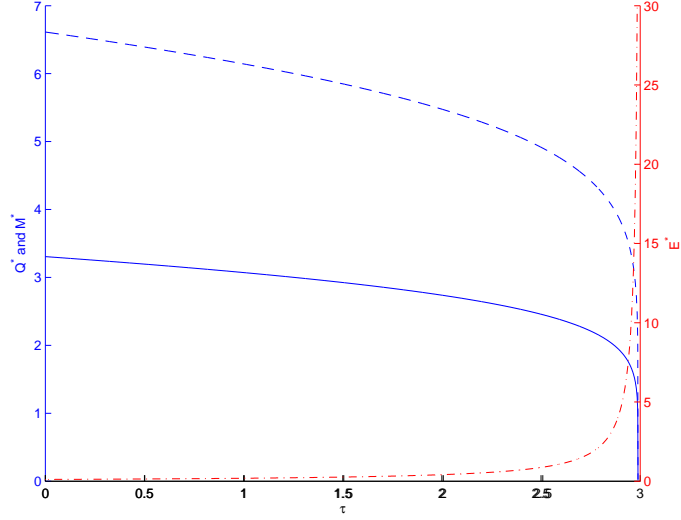}
\end{center}
\caption{The steady-states $Q^*$ (solid line), $M^*$ (dashed line), and $E^*$ (dotted line) of
system (\ref{system}) are drawn on the interval $[0,\tau_{max})$, with $\tau_{max}=2.99$, for the
parameters values given in (\ref{val1}), (\ref{val2}) and (\ref{val3}). When $\tau=\tau_{max}$,
$E^*\approx 2346$ but we have stopped the scale on the vertical axis at 30 to improve the
illustration clarity.}\label{fig_ee}
\end{figure}

With the above choices for the functions $\beta$, $g$ and $f$, we can explicitly compute the
steady-states of system (\ref{system}), $Q^*$, $M^*$ and $E^*$. In particularly, one can check that
condition (\ref{hypbeta}) holds true. Condition (\ref{condexist2}) becomes
\begin{displaymath}
(\delta+G)(a+k)<\beta_0a \quad \textrm{and} \quad
0\leq\tau<\frac{1}{\gamma}\ln\left(\frac{2\beta_0a}{(\delta+G)(a+k)+\beta_0a}\right):=\tau_{max}.
\end{displaymath}
We set
\begin{displaymath}
\alpha(\tau)=2e^{-\gamma\tau}-1 \qquad \textrm{ for } \tau\in[0,\tau_{max}).
\end{displaymath}
The function $\alpha$ is positive and decreasing on $[0,\tau_{max})$ and satisfies
\begin{displaymath}
\frac{(\delta+G)(a+k)}{a\beta_0}<\alpha(\tau)\leq 1 \qquad \textrm{ for } \tau\in[0,\tau_{max}).
\end{displaymath}
The steady-states of (\ref{system}) are then defined by
\begin{displaymath}
\begin{array}{rcl}
Q^*&=&\displaystyle\frac{\mu}{G}\frac{1}{K^{1/r}}\left(\frac{a\beta_0\alpha(\tau)-(\delta+G)(a+k)}{k(\delta+G)}\right)^{1/r},\vspace{1ex}\\
M^*&=&\displaystyle\frac{G}{\mu}Q^*,\vspace{1ex}\\
E^*&=&\displaystyle\frac{\delta+G}{\beta_0\alpha(\tau)-(\delta+G)}.
\end{array}
\end{displaymath}
For $\tau\in[0,\tau_{max})$, $Q^*$ and $M^*$ are decreasing with
\begin{displaymath}
0<Q^*\leq\frac{\mu}{G}\frac{1}{K^{1/r}}\left(\frac{a\beta_0-(\delta+G)(a+k)}{k(\delta+G)}\right)^{1/r}
\end{displaymath}
and
\begin{displaymath}
0<M^*\leq\frac{1}{K^{1/r}}\left(\frac{a\beta_0-(\delta+G)(a+k)}{k(\delta+G)}\right)^{1/r},
\end{displaymath}
and $E^*$ is increasing with
\begin{displaymath}
\frac{\delta+G}{\beta_0-(\delta+G)}\leq E^*<\frac{a}{k}.
\end{displaymath}
For the parameters given in (\ref{val1}), (\ref{val2}) and (\ref{val3}), the steady-states are
drawn on the interval $[0,\tau_{max})$ on Figure \ref{fig_ee}. In this case, $\tau_{max}\approx
2.99$.

The coefficients $A$, $B$, $C$ and $D$, defined in (\ref{coefAB}), become
\begin{displaymath}
\begin{array}{rclcrcl}
A&=&\delta+G+\beta(E^*)>0,&\qquad&C&=&\beta^{\prime}(E^*)Q^*>0,\vspace{1ex}\\
B&=&2e^{-\gamma\tau}\beta(E^*)>0,&\qquad&D&=&2e^{-\gamma\tau}\beta^{\prime}(E^*)Q^*>0,
\end{array}
\end{displaymath}
and are all strictly positive. The coefficient $G$ is constant and $H$ is still given by
$H=-f^{\prime}(M^*)>0$. One can also check that $E^*$ is the unique solution of
\begin{displaymath}
\left(2e^{-\gamma\tau}-1\right)\beta\left(E^*\right)=\delta+G.
\end{displaymath}
Thus,
\begin{displaymath}
A=B=(\delta+G)\frac{\alpha(\tau)+1}{\alpha(\tau)}.
\end{displaymath}
In particularly, we deduce that
\begin{displaymath}
\begin{array}{rcl}
b_1(\tau)&=&\mu^2+k^2>0,\vspace{1ex}\\
b_2(\tau)&=&\mu^2k^2+2G H\left[C(\mu+k+A)-AD\right],\vspace{1ex}\\
b_3(\tau)&=&G H(D-C)\left(2\mu kA-G H(C+D) \right).
\end{array}
\end{displaymath}
One can notice that $b_1$ is now independent of the delay $\tau$. Moreover, since $b_1>0$, the
polynomial function $h$, defined in (\ref{h}), has strictly positive roots if and only if (see
Lemma \ref{lemerp} and Lemma \ref{b3pos}) $b_3(\tau)<0$ or $b_3(\tau)\geq0$, $b_2(\tau)<0$ and
\begin{displaymath}
2\Delta(\tau)z_0(\tau)+b_1(\tau)b_2(\tau)-9b_3(\tau)>0.
\end{displaymath}

\begin{figure}[pt]
\begin{center}
\includegraphics[width=6cm, height=4cm]{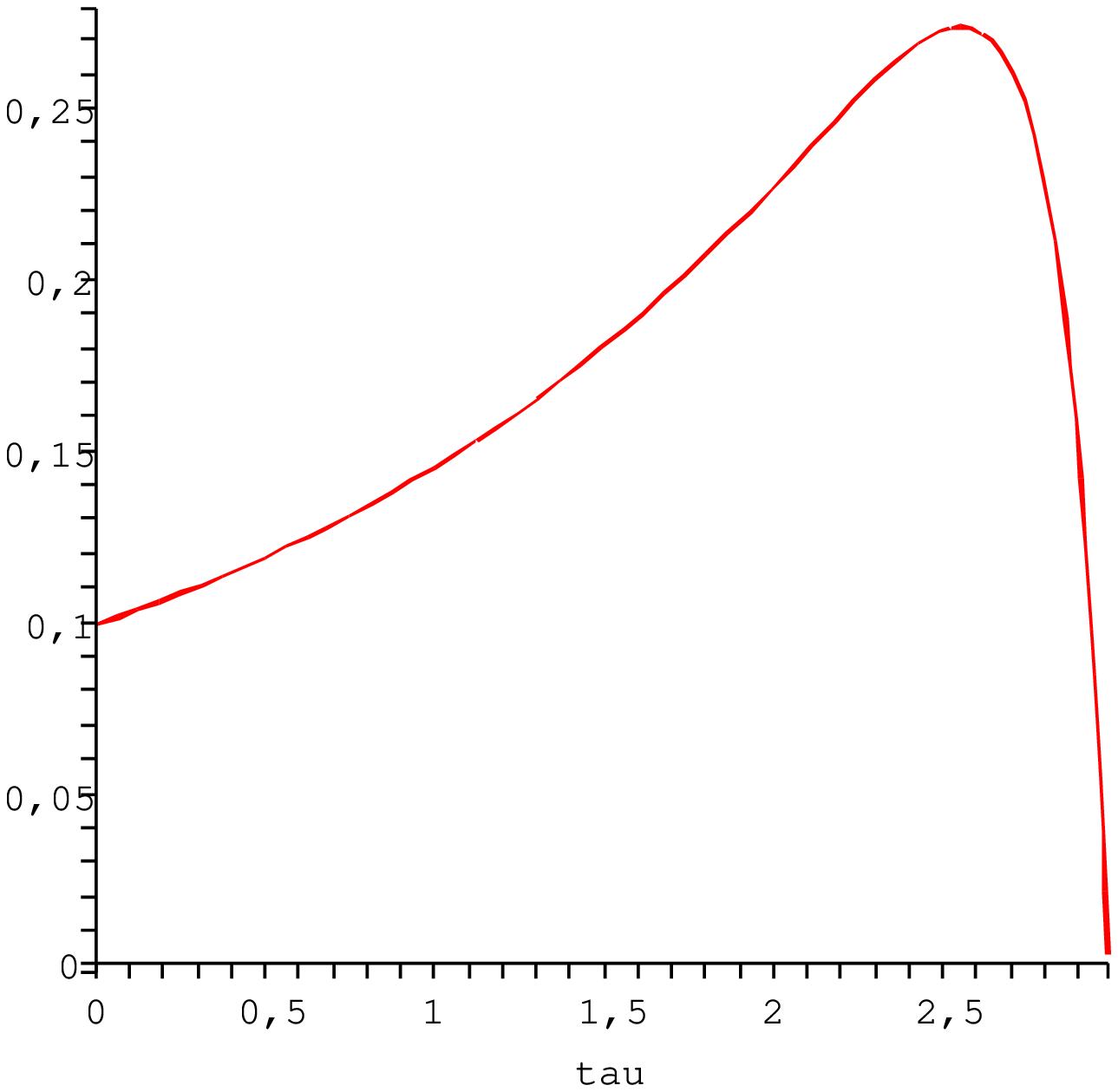}
\includegraphics[width=6cm, height=4cm]{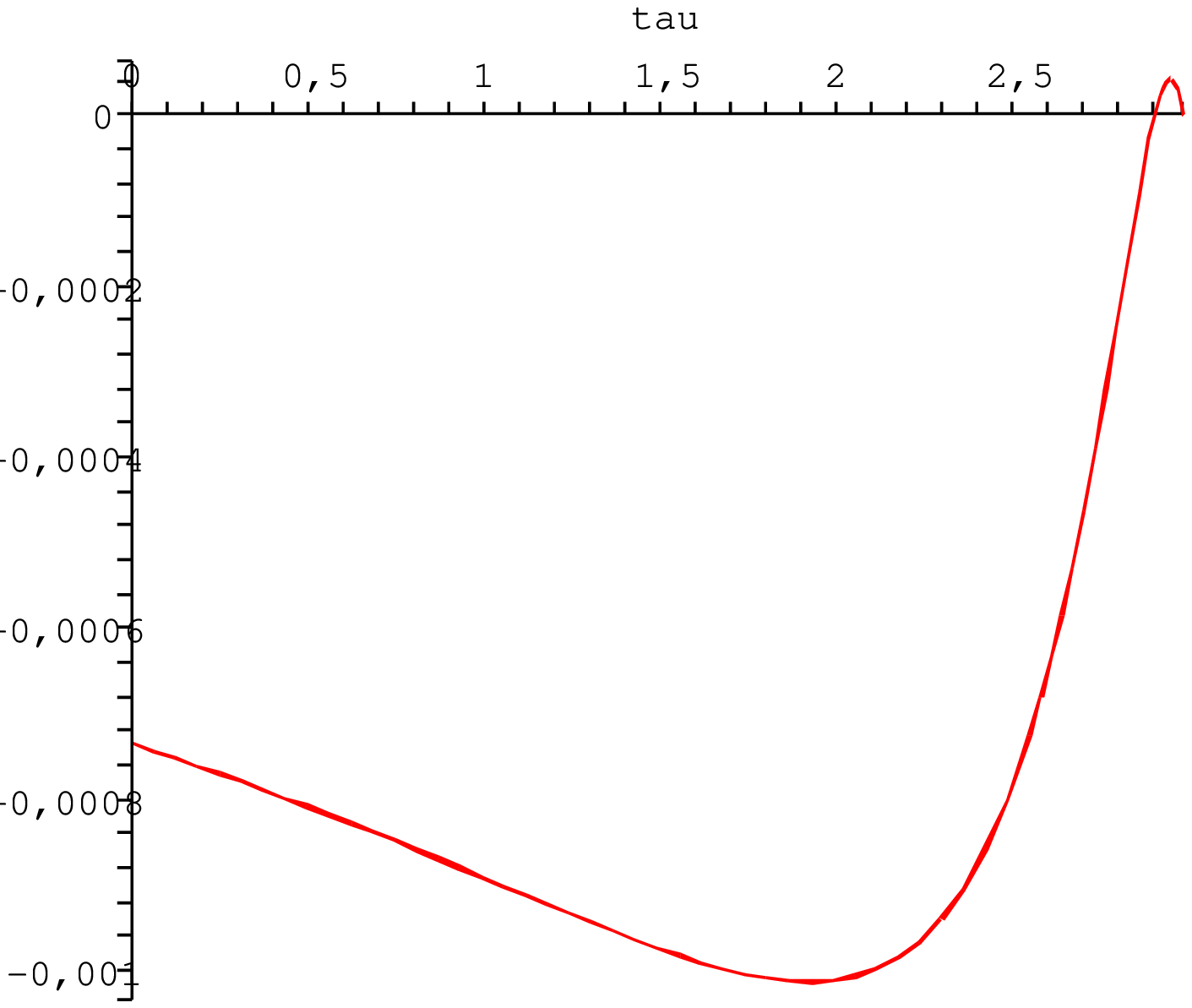}
\end{center}
\caption{Coefficients $b_2(\tau)$ (left) and $b_3(\tau)$ (right) are represented for
$\tau\in[0,\tau_{max})$ with $\tau_{max}=2.99$.}\label{fig_b}
\end{figure}

Using Maple 9, we compute the coefficients $b_2$ and $b_3$ for the values in (\ref{val1}),
(\ref{val2}) and (\ref{val3}). Results are presented in Figure \ref{fig_b}.
Since $b_3<0$ on $[0,2.92)$ and $b_2$ is always positive, $h$ has positive roots if and only if
$\tau\in I := [0,2.92)$. In this case, $h$ has exactly one positive root for $\tau^*\in[0,2.92)$,
denoted by $z^*$, and, since $h(0,\tau)<0$, $z^*$ satisfies
\begin{displaymath}
\frac{\partial h}{\partial z}(z^*,\tau^*)>0.
\end{displaymath}

\begin{figure}[hpt]
\begin{center}
\includegraphics[width=9cm, height=6cm]{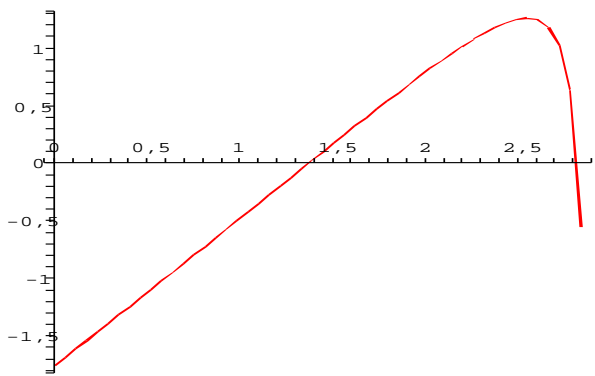}
\end{center}
\caption{Graph of the function $S_0(\tau)$ for $\tau\in[0,\tau_{max})$ with parameters given by
(\ref{val1}), (\ref{val2}) and (\ref{val3}), and $\tau_{max}\approx 2.99$. Two critical values of
$\tau$, for which stability switches can occur, appear.}\label{fig1}
\end{figure}

The function $S_0$ is drawn for $\tau\in I=[0,2.92)$ in Figure \ref{fig1}. One can see that there
are two critical values of the delay $\tau$ for which stability switches occur. In particularly,
from Theorem \ref{theorembif}, a Hopf bifurcation occurs when $\tau$ is approximately equal to
$1.4$. Thus, periodic solutions appear.

\begin{figure}[!hpt]
\begin{center}
\includegraphics[width=9cm, height=6cm]{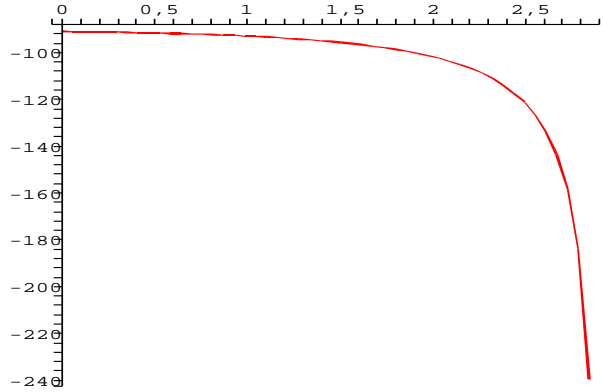}
\end{center}
\caption{Graph of the function $S_1(\tau)$ for the same values as in Figure \ref{fig1}. The
function has no positive root.}\label{fig2}
\end{figure}

In Figure \ref{fig2}, one can check that $S_1$ has no positive root on $I$. Therefore, there exist
only two critical values of the delay for which stability switches occur, $\tau=1.4$ and
$\tau=2.82$.

\begin{figure}[hpt]
\begin{center}
\includegraphics[width=9cm, height=6cm]{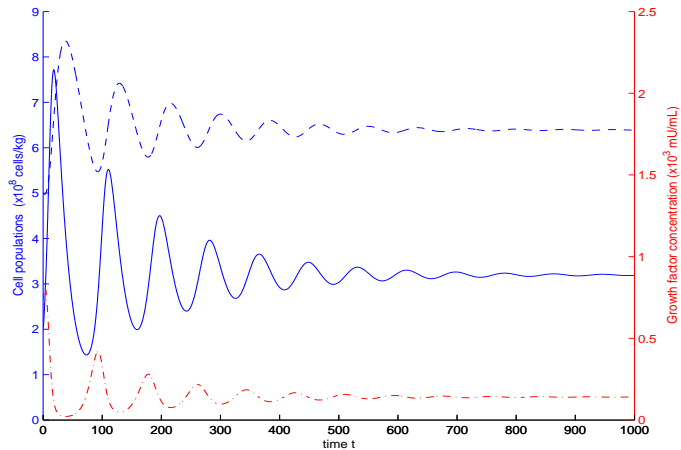}
\end{center}
\caption{Solutions $Q(t)$ (solid line), $M(t)$ (dashed line) and $E(t)$ (dotted line) of
(\ref{system}) are asymptotically stable and converge to the steady-state values. Damped
oscillations can be observed. Parameters values are the same as in Figure \ref{fig1} with
$\tau=0.5$.}\label{figstab}
\end{figure}

\begin{figure}[hpt]
\begin{center}
\includegraphics[width=9cm, height=6cm]{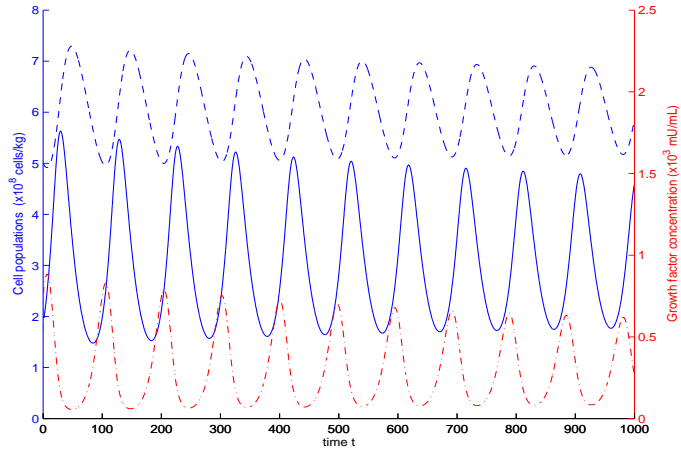}
\end{center}
\caption{When $\tau=1.4$, a Hopf bifurcation occurs and periodic solutions appear, with the same
period for the three solutions $Q(t)$ (solid line), $M(t)$ (dashed line) and $E(t)$ (dotted line)
of (\ref{system}). Periods are about 100 days. Parameters values are the same as in Figure
\ref{fig1}.}\label{fighopbif}
\end{figure}

\begin{figure}[!hpt]
\begin{center}
\includegraphics[width=9cm, height=6cm]{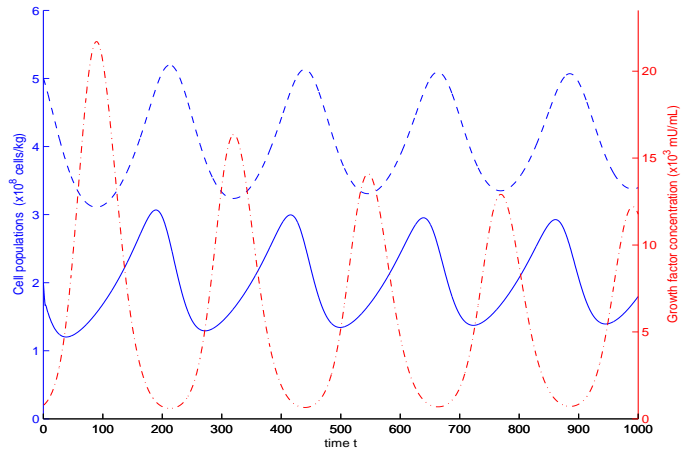}
\end{center}
\caption{For $\tau=2.8$, long period oscillations are observed, with periods close to 220 days.
Solutions $Q(t)$ (solid line), $M(t)$ (dashed line) and $E(t)$ (dotted line) of (\ref{system}) are
unstable. Parameters values are the same as in Figure \ref{fig1}.}\label{figbif2}
\end{figure}

\begin{figure}[!hpt]
\begin{center}
\includegraphics[width=9cm, height=6cm]{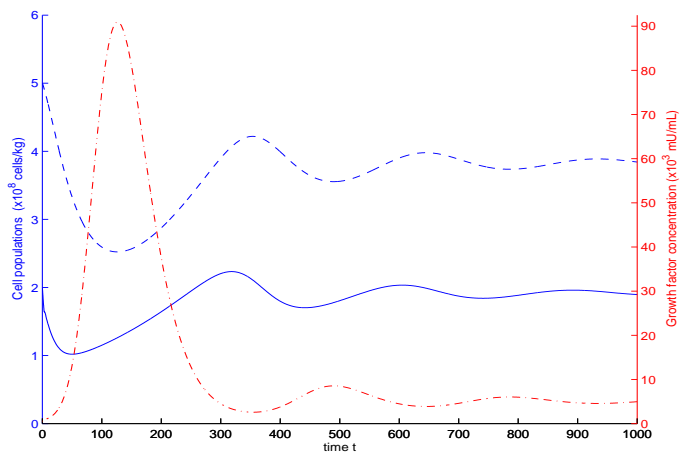}
\end{center}
\caption{When $\tau=2.9$, the steady-state is asymptotically stable and solutions $Q(t)$ (solid
line), $M(t)$ (dashed line) and $E(t)$ (dotted line) of (\ref{system}) converge to the equilibrium.
Parameters values are the same as in Figure \ref{fig1}.}\label{figstab2}
\end{figure}

Using dde23 \cite{dde23}, a {\sc Matlab} solver for delay differential equation, we can compute the
solutions of (\ref{system}) for the above mentioned values of the parameters. Illustrations are
showed in Figures \ref{figstab} to \ref{figstab2}.

Before the Hopf bifurcation occurs, solutions are stable and converge to the equilibrium, although
they oscillate transiently (see Figure \ref{figstab}). When the bifurcation occurs, periodic
solutions appear with periods about 100 days (see Figure \ref{fighopbif}). These are very long
periods compared to the delay $\tau$, which is equal to $1.4$ day.

When $\tau$ increases, periods of the oscillations increase and just before a stability switch, for
$\tau=2.82$, periods of the oscillations are close to 220 days (see Figure \ref{figbif2}). Then,
the steady-state becomes asymptotically stable again and solutions converge to the equilibrium (see
Figure \ref{figstab2}).

\section{Periodic Hematological Diseases}
Periodic hematological diseases \cite{hdm1998} represent one kind of diseases affecting blood
cells. They are characterized by significant oscillations in the number of circulating cells, with
periods ranging from weeks (19--21 days for cyclical neutropenia \cite{hdm1998}) to months (30 to
100 days for chronic myelogenous leukemia \cite{hdm1998}) and amplitudes varying from normal to low
levels or normal to high levels, depending on cell types. Because of their dynamic character,
periodic hematological diseases offer an opportunity to understand some of the regulating processes
involved in the production of hematopoietic cells.

Some periodic hematological diseases involve only one type of blood cells, for example, red blood
cells in periodic autoimmune hemolytic anemia \cite{bmm1995} or platelets in cyclical
thrombocytopenia \cite{sbmm2000}. In these cases, periods of the oscillations are usually between
two and four times the cell cycle duration. However, other periodic hematological diseases, such as
cyclical neutropenia \cite{hdm1998} or chronic myelogenous leukemia \cite{fm1999}, show
oscillations in all of the circulating blood cells, i.e., white blood cells, red blood cells and
platelets. These diseases involve oscillations with quite long periods (on the order of weeks to
months). A destabilization of the pluripotential stem cell population induced by growth factors
seems to be at the origin of these diseases.

Recently, Pujo-Menjouet {\it et al.} \cite{pm2004, pbm} and Adimy {\it et al.} \cite{acr1, acr2}
considered models for the regulation of stem cell dynamics and noticed that long periods
oscillations could be observed in hematopoiesis models. In this work, we have been able to obtain
long period oscillations (about 100 days) for a blood production model mediated by growth factors,
as a result of the feedback from blood to growth factors. This result stresses the role of growth
factors in the appearance of periodic solutions in hematopoiesis models.


One can notice that by assuming that the cell cycle duration is constant, values of $\tau$ for
which a nontrivial steady-state exists are limited and cannot be too large. This does not appear in
a model with distributed delay, as studied by Adimy {\it et al.} \cite{acr1, acr2}.

Numerical simulations demonstrated that long period oscillations in the circulating cells are
possible in our model even with short cell cycle durations. Thus, we are able to characterize some
hematological diseases, especially those exhibit a periodic behavior of all the circulating blood
cells.


\end{document}